# Factorization of the Determinant of the Gaussian-Correlation Matrix of Evenly Spaced Points Using an Inter-dimensional Multiset Duality

Selden Crary

*Palo Alto, USA*

**Abstract**

We prove our earlier conjecture that the determinant of a Gaussian-correlation matrix $\boldsymbol{V}$ with elements $V_{i.j} = \sigma_z^2 e^{-\theta(i-j)^2\delta^2}$ $(\sigma_z, \theta \in \mathbb{R};\ \sigma_z, \theta > 0;\ i, j \in \mathbb{N}_{\geq 1};\ 1 \leq i, j \leq n)$, of evenly spaced points with nearest-neighbor distance $\delta > 0$, is $\sigma_z^2 C(n)\delta^{n(n-1)}$ + higher-degree terms in $\delta$. We show that $C(n) = SF(n-1) \cdot (2\theta)^{n(n-1)/2}$, where $SF$ is the superfactorial operator. The proof uses Neville elimination to determine all elements of the upper triangular matrix $\boldsymbol{U}$ of $\boldsymbol{V}$ and provides a factorization of $det(\boldsymbol{V})$. The proof makes use of an inter-dimensional multiset duality, which involves simplices that emerge during the factorization.

## 1 Motivation

We presented a conjecture, two years ago, that the determinant of the Gaussian-correlation matrix $\boldsymbol{V}$ with elements $V_{i.j} = \sigma_z^2 e^{-\theta(i-j)^2\delta^2}$ $(\sigma_z^2, \theta \in \mathbb{R};\ \sigma_z^2, \theta > 0;\ i, j \in \mathbb{N}_{\geq 1};\ 1 \leq i, j \leq n)$, of evenly spaced points with nearest-neighbor distance $\delta > 0$, is $C(n)\delta^{n(n-1)}$ + higher-degree terms in $\delta$ [1]. This paper provides a proof of this conjecture.

## 2 Outline

After a few simple definitions and two elementary algebraic identities, we prove six multiset identities, the ultimate of which relates the union of a pair of multisets, defined on the lattice points of an $n$-simplex, to the union of a corresponding pair of multisets defined on an $(n+1)$-simplex. This identity, which can be considered variously as a lifting transformation or a duality, is key to the proof of our earlier conjecture, which we state as a lemma. The proof proceeds via Neville elimination and provides a complete determination of the upper triangular matrix $\boldsymbol{U}$ of $\boldsymbol{V}$, as well as a complete factorization of $det(\boldsymbol{V})$. We end by conjecturing that the 1D Gaussian-correlation matrix of evenly spaced points, under the conditions used in this paper, is strictly totally positive [2].

## 3 Definitions and Algebraic Identities

**Definitions:**

$\mathbb{N}_{\geq 0}$ denotes the set of whole numbers $\{0,1,2,\cdots\}$.

$\mathbb{N}_{\geq 1}$ denotes the set of natural numbers $\{1,2,3,\cdots\}$.

$\mathbb{R}$ denotes the set of real numbers.

$\eta \equiv e^{-\theta\delta^2} > 0$ $(\delta^2, \theta \in \mathbb{R};$ with the restrictions $\delta^2, \theta > 0)$.

$h_q \equiv 1 - \eta^{2q} > 0$ $(q \in \mathbb{N}_{\geq 1})$.

$SF(n) \equiv \prod_{k=1}^{n} k!$ $(n \in \mathbb{N}_{\geq 1})$ is the superfactorial operator [3].

Neville-elimination: Neville elimination is pivot-free Gaussian elimination of the $nxn$ upper triangular matrix $\boldsymbol{U}$ in LU-decomposition and uses the following formula for the relevant Stage $s+1$, Row $i$, and Column $j$ elements, in terms of elements at the immediately prior stage*:*

$U(s+1,i,j) = U(s,i,j) - \frac{U(s,i,s)U(s,s,j)}{U(s,s,s)}$ $(s,i,j \in \mathbb{N}_{\geq 1}, s \leq i,j \leq n)$ [4].

**Algebraic identities:**

**AI1:** $(i-n)^2 + (j-n)^2 = 2(i-n)(j-n) + (i-j)^2$ $(i,j,n \in \mathbb{N}_{\geq 1})$.

**AI2:** $h_{j-1} \sum_{k=0}^{i-2} \eta^{2k(j-1)} = 1 - \eta^{2(i-1)(j-1)}$ $(i,j \in \mathbb{N}_{\geq 1},\ i \geq 2,\ j \geq 1)$.

*Proof:* $h_{j-1} \sum_{k=0}^{i-2} \eta^{2k(j-1)} = \left[1 - \eta^{2(j-1)}\right]\left[1 + \eta^{2(j-1)} + \cdots + \eta^{2(i-2)(j-1)}\right]$

$= 1 - \eta^{2(i-1)(j-1)}$. □

## 4 Notation for Simplicial Multisets

Consider an array of points from a 2x1 rectangular lattice aligned commensurately with the axes and vertices of a Manhattan-aligned 2-simplex of size 6x3, as show in Fig. 1, below. The multiset of L1-norm distances from the simplex's lower-left vertex to each of the overlying lattice points is $\{0,2,3,4,5,6,6,7,8,9\}$, or written somewhat differently, $\left\{\begin{matrix} 0 \\ 2\ \ 3 \\ 4\ \ 5\ \ 6 \\ 6\ \ 7\ \ 8\ \ 9 \end{matrix}\right\}$. We are interested in generalizations of such multisets, as they will prove useful in the proof of the lemma.

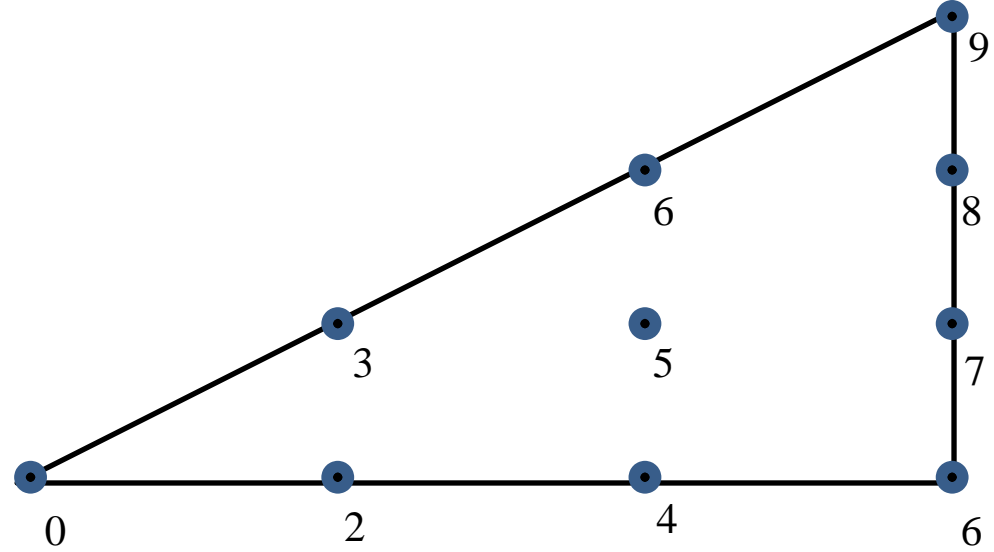

*Figure 1: A simple example of a point lattice aligned commensurately with the axes and vertices of a simplex is a rectangular lattice of points aligned with a right triangle.*

We define the following notation for the multiset of $L_1$-norm distances, added to a possibly non-zero constant $\alpha$, from an apex of an $n$-simplex to the points of a hyper-rectangular lattice aligned commensurately with the axes and vertices of the simplex, with $k, k', \cdots, k^{(n-1)'}$ being the geometric coordinates, where the notation $k^{(n-1)'}$ denotes the letter $k$ appended with $n-1$ prime symbols, except that we distinguish $k^{(n-1)'}$ from $k$:

$$\mathcal{S}_{n,\alpha,\beta,\gamma,\delta,\varepsilon,\zeta} \equiv \{\alpha + \beta k + k' + \cdots + k^{(n-1)'}\}.$$

Each unit increase of the first coordinate, viz. $k$, results in a $\beta$-unit increase in $L_1$-norm spacings. All other coordinates enter on roughly equal footing, with each unit increase in any coordinate resulting in a unit increase in $L_1$-norm spacings, all other coordinates held fixed. We assert the following conditions and ranges on the parameters.

| Conditions | Ranges | Defined/Undefined |
|---|---|---|
| $n \in \mathbb{N}_{\geq 1}$ | $k = (\gamma - 1), \cdots, (\delta - 1)$ | $k$ is always defined. |
| $\alpha, \beta \in \mathbb{N}_{\geq 0}$ | $k' = \varepsilon[k - (\varepsilon - 1)], \cdots, k - \zeta$ | $k^{(0)'}$ is undefined. |
| $\gamma, \delta \in \mathbb{N}_{\geq 1}$ | special cases for $k'$ | $k', \cdots, k^{(n-1)'}$ are undefined, if not present. |
| $\varepsilon, \zeta \in \{0,1\}$ | $(\varepsilon, \zeta) = (0,0)$: $k' = 0, \dots, k$ | |
| | " $= (1,0)$: $k' = k$ | |
| $\gamma \leq \delta$ | " $= (0,1)$: $k' = 0, \dots, k - 1$ | Undefined terms are to be neglected. |
| $\varepsilon + \zeta \in \{0,1\}$ | $k'' = 0, \cdots, k'$ | |
| | $\vdots$ | |
| | $k^{(n-1)'} = 0, \cdots, k^{(n-2)'}$ | |

Nota bena: In the above definition of $\mathcal{S}_{n,\alpha,\beta,\gamma,\delta,\varepsilon,\zeta}$, the first (resp., 2'nd, 3'rd, etc.) subscript is denoted by $n$ (resp., $\alpha, \beta, \gamma, \delta, \varepsilon, \zeta$) and must follow the conditions and ranges for $n$ (resp., $\alpha, \beta, \gamma, \delta, \varepsilon, \zeta$) given in the table, immediately above, whatever other value or symbol to which it may be assigned in any specific application. As we are interested in $\mathcal{S}$'s in which there are relationships between or among the indices, we will find notation like $\mathcal{S}_{n,\alpha,\beta,1,\delta-1,0,0}$, which simply means the fifth subscript, viz., $\delta - 1$, must follow the given rules for the fifth subscript, i.e. for what is denoted "$\delta$" in the table, above.

Detail for cases with $n = 1, 2$ or 3:

> n = 1: All $k$'s with a prime symbol are undefined, so $\mathcal{S}_{1,\alpha,\beta,\gamma,\delta,\varepsilon,\zeta} \equiv \{\alpha + \beta k\}$, with the range $(\gamma - 1) \leq k \leq (\delta - 1)$. Indices $\varepsilon$ and $\zeta$ are not used.
>
> n = 2: $\mathcal{S}_{2,\alpha,\beta,\gamma,\delta,\varepsilon,\zeta} \equiv \{\alpha + \beta k + k'\}$, with the same range as for the $n = 1$ case, as well as $\varepsilon[k - (\varepsilon - 1)] \leq k' \leq k - \zeta$.
>
> $n = 3$: $\mathcal{S}_{3,\alpha,\beta,\gamma,\delta,\varepsilon,\zeta} \equiv \{\alpha + \beta k + k' + k''\}$, with the same ranges as for $n = 2$ case, as well as $0 \leq k'' \leq k'$.

We note that if $\gamma > 1$ or $\varepsilon + \zeta \neq 0$, the relevant lattice points for $\mathcal{S}_{n,\alpha,\beta,\gamma,\delta,\varepsilon,\zeta}$ do not overlie a simplex, but just a contiguous part of one. Detail: In the former case, the lowest-valued $k$ is greater than zero, while in the latter case, the highest-valued $k'$ is less than $k$.

We make special note that one coordinate, viz., $k$, is always present and behaves differently from the rest, which appear essentially identically.

The example shown in Fig. 1, above, is, by this definition, $\mathcal{S}_{2,0,2,1,4,0,0}$. Another simple example is $\mathcal{S}_{2,0,4,1,4,0,0} = \left\{\begin{matrix} 0 \\ 4\ 5 \\ 8\ 9\ 10 \\ 12\ 13\ 14\ 15 \end{matrix}\right\}$, a multiset that will reappear, below, as an example.

$\mathcal{S}_{n,\alpha,\beta,\gamma,\delta,\varepsilon,\zeta} \uplus \mathcal{S}_{\tilde{n},\tilde{\alpha},\tilde{\beta},\tilde{\gamma},\tilde{\delta},\tilde{\varepsilon},\tilde{\zeta}}$ is the multiset union of $\mathcal{S}_{n,\alpha,\beta,\gamma,\delta,\varepsilon,\zeta}$ and $\mathcal{S}_{\tilde{n},\tilde{\alpha},\tilde{\beta},\tilde{\gamma},\tilde{\delta},\tilde{\varepsilon},\tilde{\zeta}}$.

$(-1) \otimes \mathcal{S}_{n,\alpha,\beta,\gamma,\delta,\varepsilon,\zeta}$ is the multiset $\mathcal{S}_{n,\alpha,\beta,\gamma,\delta,\varepsilon,\zeta}$ , with each element multiplied by -1.

## 5 Simplex Multiset Identities

We now give six multiset identities, named MI1 through MI6, relating $\mathcal{S}_{n,\alpha,\beta,\gamma,\delta,\varepsilon,\zeta}$ multisets.

**MI1:** $\mathcal{S}_{n,\alpha,\beta,1,\delta,0,0} = \mathcal{S}_{n,\alpha,\beta,1,\delta-1,0,0} \uplus \mathcal{S}_{n,\alpha,\beta,\delta,\delta,0,0}$ . (The RHS is the union of the two multisets resulting from slicing off the $k = \delta - 1$ face of the $n$-simplex on the LHS.)

*Example:* $(n, \alpha, \beta, \gamma, \delta, \varepsilon, \zeta) = (2,0,4,1,4,0,0)$.

$$\mathcal{S}_{n,\alpha,\beta,1,\delta,0,0} = \mathcal{S}_{2,0,4,1,4,0,0} = \left\{\begin{matrix} 0 \\ 4\ 5 \\ 8\ 9\ 10 \\ 12\ 13\ 14\ 15 \end{matrix}\right\} = \left\{\begin{matrix} 0 \\ 4\ 5 \\ 8\ 9\ 10 \end{matrix}\right\} \uplus \{12\ 13\ 14\ 15\}$$

$$= \mathcal{S}_{2,0,4,1,3,0,0} \uplus \mathcal{S}_{2,0,4,4,4,0,0} = \mathcal{S}_{n,\alpha,\beta,1,\delta-1,0,0} \uplus \mathcal{S}_{n,\alpha,\beta,\delta,\delta,0,0}.$$

$$\mathcal{S}_{n,\alpha,\beta,1,\delta,0,0} = \mathcal{S}_{2,0,4,1,4,0,0} = \left\{\begin{matrix} 0 \\ 4\ 5 \\ 8\ 9\ 10 \\ 12\ 13\ 14\ 15 \end{matrix}\right\} = \left\{\begin{matrix} 0 \\ 4\ 5 \\ 8\ 9\ 10 \end{matrix}\right\} \uplus \{12\ 13\ 14\ 15\}$$

*Proof:* By the definition of $\mathcal{S}_{n,\alpha,\beta,\gamma,\delta,0,0}$, the coordinate $k$ ranges from $0$ through $\delta - 1$. This multiset can be decomposed into the union of two mutually exclusive and exhaustive multisets, with different values of the fourth and fifth indices, and with all other indices unchanged. In the former of these constituent multisets, $k$ ranges from $0$ through $\delta - 2$, while, in the latter, $k = \delta - 1$. These two constituent multisets are $\mathcal{S}_{n,\alpha,\beta,1,\delta-1,0,0}$ and $\mathcal{S}_{n,\alpha,\beta,\delta,\delta,0,0}$, respectively. □

Three special cases that will be used later are the following:

**MI1a:** $\mathcal{S}_{n,0,\beta,1,\delta,0,0} = \mathcal{S}_{n,0,\beta,1,\delta-1,0,0} \uplus \mathcal{S}_{n,0,\beta,\delta,\delta,0,0}$ .

**MI1b:** $\mathcal{S}_{n+1,\beta-1,\beta-1,1,\delta-1,0,0} = \mathcal{S}_{n+1,\beta-1,\beta-1,1,\delta-2,0,0} \uplus \mathcal{S}_{n+1,\beta-1,\beta-1,\delta-1,\delta-1,0,0}$ .

**MI1c:** $\mathcal{S}_{n,(\beta-1)(\delta-1),1,1,\delta,0,0} = \mathcal{S}_{n,(\beta-1)(\delta-1),1,1,\delta-1,0,0} \uplus \mathcal{S}_{n,(\beta-1)(\delta-1),1,\delta,\delta,0,0}$ .

**MI2:** $\mathcal{S}_{n+1,0,\beta-1,1,\delta-1,0,0} = \mathcal{S}_{n+1,\beta-1,\beta-1,1,\delta-2,0,0} \uplus \mathcal{S}_{n+1,0,\beta-1,1,\delta-1,1,0}$ . (The RHS is the union of the two multisets resulting from slicing off all $k' < k$ of the $(n+1)$-simplex on the LHS.)

*Example:* $(n, \alpha, \beta, \gamma, \delta, \varepsilon, \zeta) = (2,0,4,1,4,0,0)$.

$$\mathcal{S}_{n+1,0,\beta-1,1,\delta-1,0,0} = \mathcal{S}_{3,0,3,1,3,0,0} = \left\{ \begin{array}{c} 0 \\ 3\;\;4\;\;5 \\ 6\;\;7\;\;8\;\;8\;\;9\;\;10 \end{array} \right\} = \left\{ \begin{array}{c} 3 \\ 6\;7\;8 \end{array} \right\} \uplus \left\{ \begin{array}{c} 0 \\ 4\;\;5 \\ 8\;\;9\;10 \end{array} \right\}$$

$$= \begin{pmatrix} \mathcal{S}_{3,3,3,1,2,0,0} \\ \uplus\, \mathcal{S}_{3,0,3,1,3,1,0} \end{pmatrix} = \begin{pmatrix} \mathcal{S}_{n+1,\beta-1,\beta-1,1,\delta-2,0,0} \\ \uplus\, \mathcal{S}_{n+1,\;0\;,\beta-1,1,\delta-1,1,0} \end{pmatrix}.$$

*Proof:* By the definition of $\mathcal{S}_{n+1,0,\beta-1,1,\delta-1,0,0}$, the coordinate $k'$ ranges from 0 through $k$. This multiset can be decomposed into the multiset union of two multisets, with different values of the sixth and seventh indices, and with all other indices unchanged. In the former of these constituent multisets, $k'$ ranges from 0 through $k-1$, while in the latter $k' = k$. These two constituent multisets are $\mathcal{S}_{n+1,0,\beta-1,1,\delta-1,0,1}$ and $\mathcal{S}_{n+1,0,\beta-1,1,\delta-1,1,0}$, respectively. However, these two multisets are not mutually exclusive, as they share the 0 element. This element can be removed from the former of the constituent multisets, if that multiset is changed to $\mathcal{S}_{n+1,\beta-1,\beta-1,1,\delta-2,0,0}$. □

**MI3:** $\mathcal{S}_{n,0,\beta,1,\delta-1,0,0} = \mathcal{S}_{n+1,0,\beta-1,1,\delta-1,1,0}$ .

*Example:* $(n, \alpha, \beta, \gamma, \delta, \varepsilon, \zeta) = (2,0,4,1,4,0,0)$. Consider the non-crossed-out integers of

$$\mathcal{S}_{n,0,\beta,1,\delta-1,0,0} = \left\{ \begin{array}{c} 0 \\ 4\;\;5 \\ 8\;\;9\;10 \\ \not{1}\!\not{2}\;\not{1}\!\not{3}\;\not{1}\!\not{4}\;\not{1}\!\not{5} \end{array} \right\} = \left\{ \begin{array}{c} 0 \\ \not{3}\;\;4\;\;5 \\ \not{6}\;\;\not{7}\;\;\not{8}\;\;8\;\;9\;\;10 \end{array} \right\} = \mathcal{S}_{n+1,0,\beta-1,1,\delta-1,1,0}.$$

*Proof:* $\mathcal{S}_{n,0,\beta,1,\delta-1,0,0} = \{\beta k + k' + \cdots + k^{(n-1)\prime}\}$, with ranges $\begin{cases} k\;\; = 0, \cdots, \delta-2 \\ k' = 0, \cdots, k \\ k'' = 0, \cdots, k' \\ \vdots \\ k^{(n-1)\prime} = 0, \cdots, k^{(n-2)\prime} \end{cases}$ . (MI3.1)

$$\mathcal{S}_{n+1,0,\beta-1,1,\delta-1,1,0} = \{(\beta-1)k + k' + \cdots + k^{n\prime}\}, \quad \begin{cases} k\;\; = 0, \cdots, \delta-2 \\ k' = k \\ k'' = 0, \cdots, k' \\ k''' = 0, \cdots, k'' \\ \vdots \\ k^{(n)\prime} = 0, \cdots, k^{(n-1)\prime} \end{cases}$$

$$= \{(\beta-1)k + k + k'' + \cdots + k^{n\prime}\}, \quad \begin{cases} k \;= 0,\cdots,\delta-2 \\ k'' = 0,\cdots,k' = k \\ k''' = 0,\cdots,k'' \\ \vdots \\ k^{(n)\prime} = 0,\cdots,k^{(n-1)\prime} \end{cases}$$

$= \{\beta k + k'' + \cdots + k^{n\prime}\}$, with the same ranges as in the last row.

Then, after a change of coordinates, in the last curly brackets, from $[k'',\cdots,k^{n\prime}]$ to $[k',k'',\cdots,k^{(n-1)\prime}]$,

$$\mathcal{S}_{n+1,\alpha,\beta-1,1,\delta-1,1,0} = \{\beta k + k' + \cdots + k^{(n-1)\prime}\}, \quad \begin{cases} k \;= 0,\cdots,\delta-2 \\ k' = 0,\cdots,k \\ k'' = 0,\cdots,k' \\ \vdots \\ k^{(n-1)\prime} = 0,\cdots,k^{(n-2)\prime} \end{cases},$$

which agrees with Eq. MI3.1. □

**MI4:** $\mathcal{S}_{n,0,\beta,\delta,\delta,0,0} = \mathcal{S}_{n,(\beta-1)(\delta-1),1,\delta,\delta,0,0}$.

*Example:* $(n,\alpha,\beta,\gamma,\delta,\varepsilon,\zeta) = (2,0,4,1,4,0,0)$. Consider the non-crossed-out integers of

$$\mathcal{S}_{n,0,\beta,\delta,\delta,0,0} = \left\{\begin{matrix} \cancel{0} \\ \cancel{4\ 5} \\ \cancel{8\ 9\ 10} \\ 12\ 13\ 14\ 15 \end{matrix}\right\} = \left\{\begin{matrix} \cancel{9} \\ \cancel{10\ 11} \\ \cancel{11\ 12\ 13} \\ 12\ 13\ 14\ 15 \end{matrix}\right\} = \mathcal{S}_{n,(\beta-1)(\delta-1),1,\delta,\delta,0,0}.$$

*Proof:* $\mathcal{S}_{n,0,\beta,\delta,\delta,0,0} = \{\beta k + k' + \cdots + k^{(n-1)\prime}\}, \quad \begin{cases} k \;= \delta-1 \\ k' = 0,\cdots,k \\ k'' = 0,\cdots,k' \\ \vdots \\ k^{(n-1)\prime} = 0,\cdots,k^{(n-2)\prime} \end{cases}$

$$= \{\beta(\delta-1) + k' + \cdots + k^{(n-1)\prime}\}, \quad \begin{cases} k' \;= 0,\cdots\delta-1 \\ k'' = 0,\cdots,k \\ k''' = 0,\cdots,k'' \\ \vdots \\ k^{(n-1)\prime} = 0,\cdots,k^{(n-2)\prime} \end{cases}. \qquad \text{(MI4.1)}$$

$\mathcal{S}_{n,(\beta-1)(\delta-1),1,\delta,\delta,0,0}$

$$= \left\{(\beta-1)(\delta-1)+k+k'+\cdots+k^{(n-1)\prime}\right\}, \quad \begin{cases} k = \delta-1 \\ k' = 0,\cdots,k \\ k'' = 0,\cdots,k' \\ \vdots \\ k^{(n-1)\prime} = 0,\cdots,k^{(n-2)\prime} \end{cases}$$

$$= \left\{(\beta-1)(\delta-1)+(\delta-1)+k'+\cdots+k^{(n-1)\prime}\right\}, \quad \begin{cases} k' = 0,\cdots,\delta-1 \\ k'' = 0,\cdots,k' \\ k''' = 0,\cdots,k'' \\ \vdots \\ k^{(n-1)\prime} = 0,\cdots,k^{(n-2)\prime} \end{cases}$$

$= \left\{\beta(\delta-1)+k'+\cdots+k^{(n-1)\prime}\right\}$, with the same ranges as in the last row.

This agrees with Eq. MI4.1. □

**MI5:** $\mathcal{S}_{n+1,\beta-1,\beta-1,\delta-1,\delta-1,0,0} = \mathcal{S}_{n,(\beta-1)\cdot(\delta-1),1,1,\delta-1,0,0}$ , $\delta \geq 2$ so that $fourth_index - 1 \geq 0$.

*Example:* $(n,\alpha,\beta,\gamma,\delta,\varepsilon,\zeta) = (2,0,4,1,4,0,0)$. Consider the non-crossed-out integers of

$$\mathcal{S}_{n+1,\beta-1,\beta-1,\delta-1,\delta-1,0,0} = \begin{Bmatrix} \cancel{3} \\ \cancel{6\ 7\ 8} \\ 9\ 10\ 11\ 11\ 12\ 13 \end{Bmatrix} = \begin{Bmatrix} 9 \\ 10\ 11 \\ 11\ 12\ 13 \\ \cancel{12\ 13\ 14\ 15} \end{Bmatrix} = \mathcal{S}_{n,(\beta-1)\cdot(\delta-1),1,1,\delta-1,0,0}.$$

*Proof:*

$$\mathcal{S}_{n+1,\beta-1,\beta-1,\delta-1,\delta-1,0,0} = \{\beta-1+(\beta-1)k+k'+\cdots+k^{n\prime}\}, \quad \begin{cases} k = \delta-2 \\ k' = 0,\cdots,k \\ k'' = 0,\cdots,k' \\ \vdots \\ k^{(n)\prime} = 0,\cdots,k^{(n-1)\prime} \end{cases}$$

$$= \{\beta-1+(\beta-1)(\delta-2)+k'+\cdots+k^{n\prime}\}, \quad \begin{cases} k' = 0,\cdots,\delta-2 \\ k'' = 0,\cdots,k' \\ k''' = 0,\cdots,k'' \\ \vdots \\ k^{(n)\prime} = 0,\cdots,k^{(n-1)\prime} \end{cases}$$

$= \{(\beta-1)(\delta-1)+k'+\cdots+k^{n\prime}\}$ , with the same ranges as in the last row. (MI5.1)

$$\mathcal{S}_{n,(\beta-1)(\delta-1),1,1,\delta-1,0,0} = \left\{(\beta-1)(\delta-1)+k+k'+\cdots+k^{(n-1)'}\right\}, \left\{\begin{array}{r}k \;=0,\cdots,\delta-2\\ k'=0,\cdots,k\\ k''=0,\cdots,k'\\ \vdots\\ k^{(n-1)'}=0,\cdots,k^{(n-2)'}\end{array}\right. .$$

Then, after changing coordinates, in the ultimate curly brackets, from $\left[k,k',\cdots,k^{(n-1)'}\right]$ to $\left[k',k'',\cdots,k^{(n)'}\right]$,

$$\mathcal{S}_{n,\beta(\delta+1),1,1,\delta-1,0,0} = \{(\beta-1)(\delta-1)+k'+\cdots+k^{n'}\}, \left\{\begin{array}{r}k'=0,\cdots,\delta-2\\ k''=0,\cdots,k'\\ k'''=0,\cdots,k''\\ \vdots\\ k^{(n)'}=0,\cdots,k^{(n-1)'}\end{array}\right. ,$$

which is the same as Eq. MI5.1. □

**MI6:** (Inter-dimensional simplex duality):

$\mathcal{S}_{n,0,\beta,1,\delta,0,0} \uplus \mathcal{S}_{n+1,\beta-1,\beta-1,1,\delta-1,0,0} = \mathcal{S}_{n+1,0,\beta-1,1,\delta-1,0,0} \uplus \mathcal{S}_{n,(\beta-1)(\delta-1),1,1,\delta,0,0}$.

*Example:* $(n,\alpha,\beta,\gamma,\delta,\varepsilon,\zeta)=(2,0,4,1,4,0,0)$.

$$\left\{\begin{array}{c}0\\ 4\;5\\ 8\;9\;10\\ 12\;13\;14\;15\end{array}\right\} \uplus \left\{\begin{array}{c}3\\ 6\;7\;8\\ 9\;10\;11\;11\;12\;13\end{array}\right\} = \left\{\begin{array}{c}0\\ 3\;4\;5\\ 6\;7\;8\;8\;9\;10\end{array}\right\} \uplus \left\{\begin{array}{c}9\\ 10\;11\\ 11\;12\;13\\ 12\;13\;14\;15\end{array}\right\}.$$

*Proof:* In sequence from left to right, the four simplex multisets in the statement can be expanded, using MI1a, MI1b, MI2, and MI1c, respectively, and the statement can be rewritten as

$$\left(\begin{array}{l}\quad\mathcal{S}_{n,0,\beta,1,\delta-1,0,0}\\ \uplus\,\mathcal{S}_{n,0,\beta,\delta,\delta,0,0}\\ \uplus\,\mathcal{S}_{n+1,\beta-1,\beta-1,1,\delta-2,0,0}\\ \uplus\,\mathcal{S}_{n+1,\beta-1,\beta-1,\delta-1,\delta-1,0,0}\end{array}\right) = \left(\begin{array}{l}\quad\mathcal{S}_{n+1,\beta-1,\beta-1,1,\delta-2,0,0}\\ \uplus\,\mathcal{S}_{n+1,0,\beta-1,1,\delta-1,1,0}\\ \uplus\,\mathcal{S}_{n,(\beta-1)(\delta-1),1,1,\delta-1,0,0}\\ \uplus\,\mathcal{S}_{n,(\beta-1)(\delta-1),1,\delta,\delta,0,0}\end{array}\right). \tag{1}$$

In the $(n,\alpha,\beta,\gamma,\delta,\varepsilon,\zeta)=(2,0,4,1,4,0,0)$ example, the first four terms of each side of Eq. 1 are colored respectively red, orange, green, and blue, in the following:

$$\left\{\begin{array}{c}0\\ 4\;5\\ 8\;9\;10\\ 12\;13\;14\;15\end{array}\right\} \uplus \left\{\begin{array}{c}3\\ 6\;7\;8\\ 9\;10\;11\;11\;12\;13\end{array}\right\} = \left\{\begin{array}{c}0\\ 3\;4\;5\\ 6\;7\;8\;8\;9\;10\end{array}\right\} \uplus \left\{\begin{array}{c}9\\ 10\;11\\ 11\;12\;13\\ 12\;13\;14\;15\end{array}\right\}. \tag{2}$$

The third multiset of the LHS of Eq. 1 cancels with the first multiset of the RHS, so Eq. 1 can be rewritten as

$$\begin{pmatrix} \mathcal{S}_{n,0,\beta,1,\delta-1,0,0} \\ \uplus\ \mathcal{S}_{n,0,\beta,\delta,\delta,0,0} \\ \uplus\ \mathcal{S}_{n+1,\beta-1,\beta-1,\delta-1,\delta-1,0,0} \end{pmatrix} = \begin{pmatrix} \mathcal{S}_{n+1,0,\beta-1,1,\delta-1,1,0} \\ \uplus\ \mathcal{S}_{n,(\beta-1)(\delta-1),1,1,\delta-1,0,0} \\ \uplus\ \mathcal{S}_{n,(\beta-1)(\delta-1),1,\delta,\delta,0,0} \end{pmatrix}.$$

In the example, this represents the cancellation of the green terms in Eq. 2.

Then, MI3, MI4, MI5 can be applied, in turn, to show that each multiset on the LHS has exactly one equal multiset on the RHS, and vice versa, as follows: the first multisets on the LHS and RHS (red, in the example) are equal, the second multiset on the LHS is equal to the ultimate multiset on the RHS (yellow, in the example), and the ultimate multiset on the LHS is equal to the second multiset on the RHS (blue, in the example). □

# 6 Lemma (Factorization of the Determinant of 1D Gaussian-Covarlance Matrices)

*a. At the end of s stages of Neville elimination, the elements of the upper triangular matrix* $\boldsymbol{U}$*, for a Gaussian-correlation matrix* $\boldsymbol{V}$ *with elements* $V_{i.j} = \sigma_z^2 e^{-\theta(i-j)^2\delta^2}$ $(\sigma_z, \theta \in \mathbb{R};\ \sigma_z, \theta > 0;\ i, j \in \mathbb{N}_{\geq 1};\ 1 \leq i, j \leq n)$ *of n evenly spaced points with nearest-neighbor distance* $\delta > 0$*, are the following, where i and j are the integer row and column indices of* $\boldsymbol{U}$*, respectively:*

| $Stage\ s$ | $Range\ of\ row\ i$ | $Range\ of\ col.\ j$ | $Element\ U(s,i,j)/\sigma_z^2$ |
|---|---|---|---|
| 1 | $1 \leq i \leq n$ | $1 \leq j \leq n$ | $\eta^{(i-j)^2}$ |
| 2 | $1 \leq i \leq 1$ | $1 \leq j \leq n$ | *same as Stage 1* |
| | $2 \leq i \leq n$ | $1 \leq j \leq 1$ | 0 |
| | | $2 \leq j \leq n$ | $\eta^{(i-j)^2} h_{j-1} \sum_{k=0}^{i-2} \eta^{2[k(j-2)+k]}$ |
| 3 | $1 \leq i \leq 2$ | $1 \leq j \leq n$ | *same as Stage 2* |
| | $3 \leq i \leq n$ | $1 \leq j \leq 2$ | 0 |
| | | $3 \leq j \leq n$ | $\eta^{(i-j)^2} h_{j-2} h_{j-1} \sum_{k=0}^{i-3} \sum_{k'=0}^{k} \eta^{2[k(j-3)+k+k']}$ |
| ⋮ | ⋮ | ⋮ | ⋮ |
| $s$ | $1 \leq i \leq s-1$ | $1 \leq j \leq n$ | *same as Stage s-1* |
| | $s \leq i \leq n$ | $1 \leq j \leq s-1$ | 0 |
| | | $s \leq j \leq n$ | $\eta^{(i-j)^2} \prod_{x=1}^{s-1} h_{j-x} \sum_{k=0}^{i-s} \sum_{k'=0}^{k} \cdots$ |
| | | | $\cdots \sum_{k^{(s-2)'}}^{k^{(s-3)'}} \eta^{2[k(j-s)+k+k'+\cdots+k^{(s-2)'}]}$. |

*b. The determinant of* $\boldsymbol{V}/\sigma_z^2$ *can be factored as* $det(\boldsymbol{V}/\sigma_z^2) = \prod_{s=2}^{n} \prod_{x=1}^{s-1} h_x$ *,which is positive.*

*c. The lowest-degree term in the expansion of* $det(\boldsymbol{V}/\sigma_z^2)$*, in powers of* $\delta$*, is* $det(\boldsymbol{V}/\sigma_z^2) = SF(n-1) \cdot (2\theta)^{n(n-1)/2} \cdot \delta^{n(n-1)}$*, where SF is the superfactorial operator.*

*d.* $\boldsymbol{V}$ *is strictly totally positive (STP), i.e., each square submatrix of* $\boldsymbol{V}$ *that consists of consecutive rows and columns has a positive determinant.*

*Proof of Part a:* The elements of $\boldsymbol{V}$ are $V_{i,j} \equiv \sigma_z^2 e^{-\theta(i-j)^2\delta^2} = \sigma_z^2 \eta^{(i-j)^2}$, where $\eta \equiv e^{-\theta\delta^2}$. Neville elimination is carried out in $n$ stages, with the first stage just copying the elements of $\boldsymbol{V}$. At the end of Stage 1, $U(1,i,j)/\sigma_z^2 = \eta^{(i-j)^2}$, for $1 \le i,j \le n$, which is the result sought for this stage. For conciseness, and without loss of generality, we drop the factor $\sigma_z^2$ that is common to all elements, in most of the remainder of the proof of this part.

We proceed with a proof by induction. During Stage 2, the first row is unchanged, while the Neville-elimination formula is applied to all other elements. At the end of Stage 2, excluding the first row, the first column is zero, and the other elements are $U(2,i,j) = U(1,i,j) - \frac{U(1,i,1)}{U(1,1,1)} U(1,1,j) = \eta^{(i-j)^2} - \eta^{(i-1)^2+(j-1)^2}$. Applying algebraic identity AI1 of Sec. 3, with $n = 1$, gives $U(2,i,j) = \eta^{(i-j)^2}\left[1 - \eta^{2(i-1)(j-1)}\right]$. Next, applying AI2 gives $U(2,i,j) = \eta^{(i-j)^2} h_{j-1} \sum_{k=0}^{i-2} \eta^{2[k(j-2)+k]}$, which is the result sought for this stage.

To complete the proof, we show that if the statement is true for arbitrary Stage $0 \le w \le n-1$, then the statement is true for Stage $(w+1)$. For Stage $w$, the statement is

$$U(w,i,j) = \eta^{(i-j)^2} \prod_{x=1}^{w-1} h_{j-x} \sum_{k=0}^{i-w} \sum_{k'=0}^{k} \cdots \sum_{k^{(w-2)'}=0}^{k^{(w-3)'}} \eta^{2\left[k(j-w)+k+k'+\cdots+k^{(w-2)'}\right]}.$$

During Stage $(w+1)$, the first $w$ rows and $w-1$ columns are unchanged, while the Neville-elimination formula $U(w+1,i,j) = U(w,i,j) - \frac{U(w,i,w)}{U(w,w,w)} U(w,w,j)$ is applied to all other elements, giving, after cancellations,

$$\begin{aligned} U(w+1,i,j) = {} & \eta^{(i-j)^2} \prod_{x=1}^{w-1} h_{j-x} \sum_{k=0}^{i-w} \sum_{k'=0}^{k} \cdots \sum_{k^{(w-2)'}=0}^{k^{(w-3)'}} \eta^{2\left[k(j-w)+k+k'+\cdots+k^{(w-2)'}\right]} \\ & - \eta^{(i-w)^2+(j-w)^2} \prod_{x=1}^{w-1} h_{j-x} \sum_{k=0}^{i-w} \sum_{k'=0}^{k} \cdots \sum_{k^{(w-2)'}=0}^{k^{(w-3)'}} \eta^{2\left[k+k'+\cdots+k^{(w-2)'}\right]}. \end{aligned}$$

Collecting common factors, and applying AI1, with $n$ in that identity being $w$ here, the last equation becomes,

$$U(w+1,i,j) = \eta^{(i-j)^2} \prod_{x=1}^{w-1} h_{j-x} \sum_{k=0}^{i-w} \sum_{k'=0}^{k} \cdots \sum_{k^{(w-2)'}}^{k^{(w-3)'}} \begin{pmatrix} \eta^{2\left[k(j-w)+k+k'+\cdots+k^{(w-2)'}\right]} \\ -\eta^{2(i-w)(j-w)} \eta^{2\left[k+k'+\cdots+k^{(w-2)'}\right]} \end{pmatrix}$$

$$(w+1 \le i,j \le n). \quad (3)$$

Applying MI6, with $(w-1, i-w+1, j-w+1)$ here being $(n,\delta,\beta)$ in the identity, gives

$$\mathcal{S}_{w-1,0,j-w+1,1,i-w+1,0,0} \uplus \mathcal{S}_{w,j-w,j-w,1,i-w,0,0} = \mathcal{S}_{w,0,j-w,1,i-w,0,0} \uplus \mathcal{S}_{w-1,(i-w)(j-w),1,1,i-w+1,0,0},$$

or allowing for negative elements in the multisets [5], and after collecting the $(w-1)$-simplices and $w$-simplices on the LHS and RHS, respectively, we get the lifting transformation,

$$\mathcal{S}_{w-1,0,j-w+1,1,i-w+1,0,0} \uplus (-1) \otimes \mathcal{S}_{w-1,(i-w)(j-w),1,1,i-w+1,0,0} = \mathcal{S}_{w,0,j-w,1,i-w,0,0} \uplus (-1) \otimes \mathcal{S}_{w,j-w,j-w,1,i-w,0,0}.$$

The LHS of the last equation is the set of exponents of $\eta^2$ in the sum of the two terms in the parenthesis in Eq. 3. Thus, Eq. 3 can be rewritten, using the RHS of the last equation,

$$
\begin{aligned}
&U(w+1,i,j) \\
&\quad = \eta^{(i-j)^2} \prod_{x=1}^{w-1} h_{j-x} \sum_{k=0}^{i-w-1} \sum_{k'=0}^{k} \cdots \sum_{k^{(w-1)'}}^{k^{(w-2)'}} \begin{pmatrix} \eta^{2[k(j-w-1)+k+k'+\cdots+k^{(w-1)'}]} \\ -\eta^{2(j-w)} \eta^{2[k(j-w-1)+k+k'+\cdots+k^{(w-1)'}]} \end{pmatrix} \\
&\quad = \eta^{(i-j)^2} \left[1-\eta^{2(j-w)}\right] \prod_{x=1}^{w-1} h_{j-x} \sum_{k=0}^{i-w-1} \sum_{k'=0}^{k} \cdots \sum_{k^{(w-1)'}}^{k^{(w-2)'}} \eta^{2[k(j-w-1)+k+k'+\cdots+k^{(w-1)'}]} \\
&\quad = \eta^{(i-j)^2} h_{j-w} \prod_{x=1}^{w-1} h_{j-x} \sum_{k=0}^{i-w-1} \sum_{k'=0}^{k} \cdots \sum_{k^{(w-1)'}}^{k^{(w-2)'}} \eta^{2[k(j-w-1)+k+k'+\cdots+k^{(w-1)'}]} \\
&\quad = \eta^{(i-j)^2} \prod_{x=1}^{w} h_{j-x} \sum_{k=0}^{i-w-1} \sum_{k'=0}^{k} \cdots \sum_{k^{(w-1)'}}^{k^{(w-2)'}} \eta^{2[k(j-w-1)+k+k'+\cdots+k^{(w-1)'}]},
\end{aligned}
$$

$$(w+1 \le i,j \le n).$$

Changing the name of the stage coordinate from $w+1$ to $s$, and including the $\sigma_z^2$, which for conciseness was dropped earlier, gives

$U(s,i,j)/\sigma_z^2 = \eta^{(i-j)^2} \prod_{x=1}^{s-1} h_{j-x} \sum_{k=0}^{i-s} \sum_{k'=0}^{k} \cdots \sum_{k^{(s-2)'}}^{k^{(s-3)'}} \eta^{2[k(j-s)+k+k'+\cdots+k^{(s-2)'}]}$, which is the expression in the statement. □

*Proof of Part b:* The determinant of a real, $n \times n$ matrix with positive elements, as is the case at hand, is the product of the main-diagonal elements of its upper triangular matrix [6], i.e., $det(\boldsymbol{V}/\sigma_z^2) = \prod_{s=1}^{n} U(s,s,s)/\sigma_z^2$, where we have chosen to write the product in terms of the diagonal elements established at the end of each Neville-elimination stage. From the statement of Part a, $U(1,1,1)/\sigma_z^2 = 1$, and $U(s,s,s)/\sigma_z^2 = \prod_{x=1}^{s-1} h_{s-x}$ $(s \ge 2)$, so $det(\boldsymbol{V}/\sigma_z^2) = \prod_{s=2}^{n} \prod_{x=1}^{s-1} h_{s-x}$, which is positive because each $h_x$ is positive. □

*Proof of Part c:* Substituting the polynomials $h_x$ $(x \in \mathbb{N}_{\ge 1})$ into the statement of Part b, gives

$det(\boldsymbol{V}/\sigma_z^2) = \prod_{s=2}^{n} \prod_{x=1}^{s-1} \left[1-\eta^{2(s-x)}\right]$.

By Definition 3, $\eta \equiv e^{-\theta\delta^2}$, so $1-\eta^{2(s-x)} = 1-e^{-2(s-x)\theta\delta^2}$, which, for sufficiently small $\delta$, expands to $2(s-x)\theta\delta^2 + O(\theta^2\delta^4)$, giving,

$$
\begin{aligned}
det(\boldsymbol{V}/\sigma_z^2) &= \prod_{s=2}^{n} \prod_{x=1}^{s-1} [2(s-x)\theta\delta^2 + O(\theta^2\delta^4)] \\
&= \prod_{s=2}^{n} \{(s-1)!\,(2\theta\delta^2)^{s-1} + O[(\theta\delta^2)^s]\} \\
&= SF(n-1) \cdot (2\theta\delta^2)^{n(n-1)/2} + higher\text{-}\ degree\ terms\ in\ \delta,
\end{aligned}
$$

thus proving the conjecture in [1]. □

*Proof of Part d:* Karlin demonstrated $\boldsymbol{V}$ is STP in Chpt. 1, Sec. 2 of [8]. We note that STP is elsewhere called total positivity [9].

# 7 Interpretation

Part b of the lemma tells us that the determinant is given by a product of $h$ functions, as follows. When there is just one point, i.e., for $n = 1$, $Det(\boldsymbol{V}/\sigma_Z^2) = 1$. When a second point is added, $Det(\boldsymbol{V}/\sigma_Z^2)$ increases by a factor $h_1$ that we can consider as due to the fact that we now have a pair of points. When a third evenly spaced point is added to the previous two, $Det(\boldsymbol{V}/\sigma_Z^2)$ increases by a factor $h_1h_2$. This can be interpreted as a factor $h_1$ due to the newly created, near-neighbor pair, as well as a factor $h_2$ due to the newly created, second-nearest-neighbor pair. In the case of five evenly spaced points in 1D, with nearest-neighbor spacing $\delta$, then $Det(\boldsymbol{V}/\sigma_Z^2) = h_1^4h_2^3h_3^2h_4$. This pairwise, multiplicative, particle-interaction interpretation holds nicely for any number of evenly spaced points.

Part c of the lemma tells us that the expansion of $Det(\boldsymbol{V})$, in powers of $\delta$, commences with a term proportional to $\delta^{n(n-1)}$, i.e. a factor $\delta^2$ for each of the $\binom{n}{2} = n(n-1)/2$ pairs of points. Because $Det(\boldsymbol{V})$ enters the denominator of each element of $\boldsymbol{V}^{-1}$, we can think of each pair of points in the design as contributing a factor $\delta^2$ to the denominator of each element of $\boldsymbol{V}^{-1}$. Of course, the effect of the design also affects the numerator, so the situation is more complex than just naïvely considering only the denominator. This subject will be taken up in a more-detailed, follow-up report.

## Version History

v2: After v1 appeared, Prof. Michael L. Stein of the Univ. of Chicago pointed us to a Y2000 paper by Wei-Liem Loh and Tao-Kai Lam [7] that had proved, using a different method, Part b of the lemma presented here.

v3: Prof. (ret.) John Nuttall of Western Ontario Univ. pointed out that Karlin [8] contains a simple proof that matrix $\boldsymbol{V}$, as used in this paper, is strictly totally positive, so we have dropped the conjecture to this effect, acquiesced to his request that the author list remain unchanged, added pointers to the relevant pages of [8], and modified the Acknowledgments and References, here, accordingly. In addition, the following changes were made: the word "isosceles" in Fig. 1 was removed; $\delta$ was changed to $\gamma$ in the fourth paragraph of Section 4; all occurrences of the variable $p$ were changed to $n$, in the statement of the lemma in Section 6; "Stage $w \geq 2$" was changed to "Stage $0 \leq w \leq n-1$," in the third paragraph of the proof of the lemma in Section 6; "Stage $(w+1) \leq n$" was changed to "Stage $(w+1)$" in the same paragraph; an errant right parenthesis was removed from Sec. 7; key words were added; and "Covariance" was corrected to "Correlation," in the title.

## Acknowledgments

We thank Prof. Michael L. Stein of the Univ. of Chicago (cf. Version History, v2, above); Prof. Wei-Liem Loh of the National Univ. of Singapore for pointing out that the word "covariance" in the title of Version 1 should have been "correlation," as the quantity $\sigma_Z^2$ used in this paper is taken as a constant; and Prof. (ret.) John Nuttall of Western Ontario University for pointing us to Karlin's book [8] for the proof of Part d of the lemma (cf. Version History, v3, above).